\documentclass{article}
\usepackage{amsmath,amssymb,amsthm,latexsym,a4}
\usepackage{hyperref,rotating}
\usepackage[a4paper]{geometry}

 \newtheorem{theorem}{Theorem}
\newtheorem*{theorem*}{Theorem}
\newtheorem{lemma}{Lemma}[section]
\newtheorem*{lemma*}{Lemma}
\newtheorem{corollary}[lemma]{Corollary}
\newtheorem*{corollary*}{Corollary}
\newtheorem{proposition}[lemma]{Proposition}
\newtheorem*{proposition*}{Proposition}

\newcommand{\PA}{\mathsf{PA}}
\newcommand{\EA}{\mathsf{EA}}
\newcommand{\ISi}{\mathsf{I}\Sigma}
\newcommand{\IPi}{\mathsf{I}\Pi}
\newcommand{\BSi}{\mathsf{B}\Sigma}

\newcommand{\Rfn}{\mathsf{Rfn}}
\newcommand{\RFN}{\mathsf{RFN}}
\newcommand{\True}{\mathsf{True}}

\newcommand{\Prf}{\mathsf{Prf}}
\newcommand{\Con}{\mathsf{Con}}
\newcommand{\nCon}{\mbox{$n$-$\mathsf{Con}$}}

\newcommand{\ul}{\ulcorner}
\newcommand{\ur}{\urcorner}
\newcommand{\la}{\langle}
\newcommand{\ra}{\rangle}

\newcommand{\imp}{\rightarrow}
\newcommand{\eqv}{\leftrightarrow} 
\newcommand{\gn}[1]{\ulcorner #1 \urcorner}
\renewcommand{\phi}{\varphi}
\newcommand{\nat}{\mathbb{N}}
\newcommand{\refeq}[1]{(\ref{#1})}

\begin{document}

\title{Axiomatizing provable $n$-provability\thanks{Research financed by a grant of the Russian Science Foundation, project No.~14-50-00005.}}

\author{Evgeny Kolmakov \qquad Lev D. Beklemishev \\
Steklov Mathematical Institute of Russian Academy of Sciences \\
Gubkina str. 8, Moscow, Russia \\
\texttt{kolmakov-ea@yandex.ru} \qquad \texttt{bekl@mi.ras.ru}
}


\maketitle

\begin{abstract}
A formula $\phi$ is called \emph{$n$-provable} in a formal  arithmetical theory $S$ if $\phi$ is provable in $S$ together with all true arithmetical $\Pi_{n}$-sentences taken as additional axioms. While in general the set of all $n$-provable formulas, for a fixed $n>0$, is not recursively enumerable, the set of formulas $\phi$ whose $n$-provability is provable in a given r.e.\ metatheory $T$ is r.e. This set is deductively closed and will be, in general, an extension of $S$. We prove that these theories can be naturally axiomatized in terms of progressions of iterated local reflection principles.
In particular, the set of provably 1-provable sentences of Peano arithmetic $\PA$ can be axiomatized by $\varepsilon_0$ times iterated local reflection schema over $\PA$. Our characterizations yield additional information on the proof-theoretic strength of these theories (w.r.t.\ various measures of it) and on their axiomatizability.
We also study the question of speed-up of proofs and show that in some cases a proof of  $n$-provability of a sentence can be much shorter than its proof from iterated reflection principles.

\vspace{0.4cm}

\end{abstract}

\section{Introduction}



A lot of recent work on provability logic and its applications to the analysis of systems of arithmetic involves the notion of $n$-provability. An arithmetical formula $\phi$ is called \emph{$n$-provable} in a formal theory $T$ whose language contains that of Peano arithmetic if $\phi$ is provable in $T$ together with all true arithmetical $\Pi_{n}$-sentences taken as additional axioms. Another notion of $n$-provability has also been considered in the literature, where $\phi$ is called $n$-provable if $\phi$ is provable in $T$ using first order logic and one application of the omega-rule stated for $\Sigma_{n-1}$-formulas (cf \cite{Smo77}). The two notions coincide only for $n=1$, however the former is simpler and more useful for applications, so $n$-provability will only refer to that notion in the present paper.

The set of $n$-provable formulas of a (sound) r.e.\ theory $T$ is $\Sigma_{n+1}$-complete and is expressible by an arithmetical  $\Sigma_{n+1}$-formula usually denoted $[n]_T(x)$. This formula (provably in $T$) shares the main properties of the usual G\"odel's provability predicate, that is, L\"ob's derivability conditions, therefore the standard proof of G\"odel's second incompleteness theorem also works for $n$-provability. Hence, if a theory $T$ is $n$-consistent, its $n$-consistency is not $n$-provable.

The sentence $\neg [n]_T(\gn{\bot})$ expressing the $n$-consistency of $T$ is also known to be equivalent to the uniform reflection principle for $\Sigma_n$-formulas, $\RFN_{\Sigma_n}(T)$. This connects the study of $n$-provability to the theory of transfinite recursive progressions of axiomatic systems based on iteration of reflection principles that originated in the works of Turing~\cite{Tur39} and Feferman~\cite{Fef62}.

In provability logic, the notion of $n$-provability emerged in the work of C.~Smory\'nski (see~\cite{Smo85}) who characterized GL as the provability logic of the $n$-provability predicate.  Konstantin~Ignatiev~\cite{Ign93}, following Giorgi~Japaridze~\cite{Dzh86}, axiomatized the polymodal provability logic GLP of the $n$-provability predicates taken together for all $n\in\omega$.

The notion of $n$-provability later found several interesting applications in the study of fragments of $\PA$. It was used to characterize the fragments of arithmetic defined by parameter free induction~\cite{Bek99} and to study their properties such as complexity of axiomatization, the classes of provably total computable functions, etc. The notion of $n$-provability and the modal logic GLP also played a prominent role in the approach to the ordinal analysis of systems of arithmetic based on provability algebras~\cite{Bek04,Bek05en}.

While the set of all $n$-provable formulas (for $n>0$) is in general not recursively axiomatizable, the set of formulas $\phi$ for which $[n]_S\phi$ is provable in the ordinary sense in a given (r.e.) metatheory $T$ is a recursively enumerable set. This set is deductively closed and will be, in general, an extension of $S$ (denoted $C^n_S(T)$ below). We are interested in obtaining a natural recursive axiomatization of such a theory and to better understand its properties such as how strong it is compared to $T$ and to $S$ and
whether it is finitely axiomatizable. For example, what is the theory $C^1_\PA(\PA)$?

We notice that $C^0_\PA(\PA)$ has the same set of theorems as $\PA$. However, by Parikh's theorem there is a non-provably recursive speed-up between $\PA$ and $C^0_\PA(\PA)$. The case $n=0$ is an exception, and for all $n>0$ the theories $C^n_S(T)$ are, in general, strictly stronger than $S$.

The answers that we obtain are formulated in terms of progressions of iterated local reflection principles. For example, we show that $C^1_S(\PA)$ is equivalent to $\varepsilon_0$ times iterated local reflection schema over $S$. More generally, it turns out that the number of times the local reflection schema needs to be iterated to axiomatize the theory $C^1_S(T)$ is related to the so-called $\Sigma^0_2$-ordinal of $T$ introduced in \cite{BekVis}. We also obtain more general characterizations of theories  $C^n_S(T)$ for $n>1$ in terms of iterations of relativized local reflection principles.

The paper is organized as follows. Section 2 contains a brief summary of basic notions and the notation used in this paper.
In Section 3 we introduce the theories $C^n_S(T)$ and prove several basic results related to them.
Section 4 is devoted to the axiomatization of the theory $C^1_S(\EA)$, which is the base case for our study.
In Section 5 we are dealing with the same problem for the extensions of $\EA$ and obtain our main axiomatization results for the case of $1$-provability.
We obtain the relativization of these results to the case of $n$-provability for $n > 1$ in Section 6.
In Section 7 we prove that the natural axiomatization of $C^1_\EA(\EA)$ has superexponential speed-up over 
the axiomatization $\EA + \Rfn(\EA)$ and obtain related results for provable $n$-provability.

The question studied in this paper emerged in discussions with Volodya Shavrukov and was suggested by some (unpublished) results of Mingzhong Cai~\cite{Cai15} (cf. Proposition~\ref{cai} in this paper).  We are grateful for their valuable input.


\section{Preliminaries}
In this paper we deal with first-order theories in the language of arithmetic. As a basic theory we take \emph{Elementary Arithmetic} $\EA$
(sometimes denoted as $\mathsf{I}\Delta_0(\exp)$), that is, the first-order theory formulated in the language $0, (\cdot)', +, \times$
extended by the unary function symbol $\exp$ for the exponentiation function $2^x$.
It has standard defining axioms for these symbols and the induction schema for all bounded formulas in this language.
As usual, $\varphi$ is a \emph{bounded formula} if it contains only bounded quantifiers.
The class of all bounded formulas in the language introduced above is denoted by $\Delta_0(\exp)$, and we also call such formulas \emph{elementary}.
We denote by $\EA^+$ the extension of $\EA$ by the axiom asserting the totality of the \emph{superexponential function} $2^x_y$.

If we allow induction for all arithmetical formulas the resulting theory is \emph{Peano Arithmetic} denoted by $\PA$.
The fragment of $\PA$ obtained by restricting the induction schema to $\Sigma_n$-formulas is denoted by $\ISi_n$.
By $\ISi_n^-$ we denote the theory of parameter free induction for $\Sigma_n$-formulas.
It is clear that the sequence of theories $\{\ISi_n \mid n < \omega\}$ is monotone and $\PA = \bigcup_{n < \omega} \ISi_n$.
We also consider the theories $\mathsf{B}\Gamma$ and $\mathsf{L}\Gamma$ and their parameter free versions $\mathsf{B}\Gamma^-$ and $\mathsf{L}\Gamma^-$,
where $\Gamma$ is some class of arithmetical formulas (see, e.g., \cite{KPD}).

We assume the standard arithmetization of syntax and g{\"o}delnumbering of syntactic objects.
In particular, we write $\ul \varphi \ur$ for the (numeral of the) g{\"o}delnumber of $\varphi$.
We consider \emph{r.e. extensions} of $\EA$ and assume that each theory $T$ comes with an elementary formula $\sigma_T(x)$,
defining the set of axioms of $T$ in the standard model. Using this formula one can naturally construct the formula $\Prf_T(x, y)$
representing the relation ``$x$ codes a $T$-proof of the formula with g{\"o}delnumber $y$''.
Then the \emph{standard provability predicate} for $T$ is given by $\exists x\, \Prf_T(x, y)$, and we denote this formula by $\Box_T(y)$.
We often write $\Box_T \varphi$ instead of $\Box_T(\ul \varphi \ur)$.
Also we allow quantification over sentences or formulas in a natural way,
e.g., $\forall \varphi\, I(\varphi)$ is a shorthand for
$\forall x\, (\mathrm{Sent}(x) \imp I(x))$, where $\mathrm{Sent}(x)$ is an elementary formula
defining the set of all g{\"o}delnumbers of the arithmetical sentences.

If every theorem of $T$ is a theorem of $U$ we write $T \subseteq U$.
The formalization of this statement in arithmetic is given by $\forall \psi\, (\Box_T \psi \imp \Box_U \psi)$.
By $T \equiv U$ we mean that both inclusions $T \subseteq U$ and $U \subseteq T$ hold, that is, $T$ and $U$ are \emph{deductively equivalent}.
If the above conditions hold only for the formulas in some class $\Gamma$ we write $T \subseteq_\Gamma U$ and $T \equiv_\Gamma U$, respectively.

It is known that for $n > 0$ classes $\Pi_n$ have \emph{partial truth definitions}.
Namely, for each $n > 0$ there exists an arithmetical $\Pi_n$-formula $\True_{\Pi_n}(x)$ such that
for every $\Pi_n$-formula $\varphi(x_1, \dots, x_n)$
$$
\EA \vdash \forall x_1\dots \forall x_n\:(\varphi(x_1, \dots, x_n) \eqv \True_{\Pi_n}(\ul\varphi(\overline{x}_1, \dots, \overline{x}_n) \ur)),
$$
where the bar notation $\ul \varphi(\overline{x}) \ur$ stands for the $\EA$-definable term, representing the elementary function
that maps $n$ to $\ul\varphi(\bar{n})\ur$, where $\bar{n}$ is the term $0^{'\dots'}$ with $n$ successor symbols.
Using this formula we can formalize the notion of \emph{$n$-provability}.
We denote by $[n]_T$ the provability predicate for the theory $T$ together with all true $\Pi_n$-sentences taken as additional axioms.
This predicate can be represented in arithmetic using the corresponding partial truth definition as follows:
$$
[n]_T \varphi := \exists z \left(\True_{\Pi_n}(z) \wedge \Box_T(\True_{\Pi_n}(\overline{z}) \imp \varphi)\right).
$$
It is known that $[n]_T$ satisfies \emph{L{\"o}b's derivability conditions} provably in $\EA$ (cf \cite{Smo85,Bek05en}):
\begin{enumerate}
\item If $T \vdash \varphi$, then $\EA \vdash [n]_T\varphi$.
\item $\EA \vdash [n]_T(\varphi \imp \psi) \imp ([n]_T\varphi \imp [n]_T\psi)$.
\item $\EA \vdash [n]_T\varphi \imp [n]_T[n]_T\varphi$.
\end{enumerate}
Point 3 follows from the general fact, known as provable \emph{$\Sigma_{n+1}$-completeness}:
$$
\EA \vdash \forall x_1\dots \forall x_n\:(\sigma(x_1, \dots, x_m) \imp [n]_T\sigma(\overline{x}_1, \dots, \overline{x}_m)),
$$
whenever $\sigma(x_1, \dots, x_m)$ is a $\Sigma_{n+1}$-formula.

We write $\la n \ra_T\varphi$ for $\neg[n]_T\neg\varphi$. The formula $\la n \ra_T \top$ represents \emph{$n$-consistency} of a theory $T$
and is known to be equivalent to the \emph{uniform reflection principle} for $\Sigma_n$-formulas
$$
\RFN_{\Sigma_n}(T) \colon \forall x\, \left(\Box_T\varphi(\bar{x}) \imp \varphi(x)\right), \ \varphi \in \Sigma_n.
$$
The union of all schemata $\RFN_{\Sigma_n}(T)$ for $n>0$, that is, the \emph{full uniform reflection schema}, is denoted $\RFN(T)$.

The \emph{local reflection principle} for $T$ is the following schema
$$
\Rfn(T) \colon \Box_T\varphi \imp \varphi, \text{ $\varphi$ a sentence}.
$$
If we impose the restriction $\varphi \in \Gamma$ for some class of formulas $\Gamma$, then we obtain the \emph{partial reflection principle} denoted by $\Rfn_\Gamma(T)$.
\emph{Relativized local reflection principles} are defined analogously, but with $[n]_T$ instead of $\Box_T$. For instance, $\Rfn^n_\Gamma(T)$ denotes the schema
$$
[n]_T\varphi \imp \varphi, \ \varphi \in \Gamma.
$$

We say that $(D, \prec)$ is an \emph{elementary linear ordering} if there is a pair of elementary formulas $D(x)$ and $x \prec y$ such that
$\EA$ proves that $(D, \prec)$ is a linear ordering. It is an \emph{elementary well-ordering},
if $\prec$ is a well-ordering of $\{n \in \mathbb{N} \mid \mathbb{N} \models D(n)\}$ in the standard model of arithmetic.

We consider \emph{transfinite iterations} of the reflection schemata mentioned above along an arbitrary elementary well-ordering $(D, \prec)$.
We follow the treatment of iterated reflection principles presented in~\cite{Bek99b}.
Assuming $R(T)$ is one of the reflection schemata for $T$ a \emph{transfinite  progression} $(R(T))_{\alpha \in D}$ is
defined by formalizing the following fixed-point equation:
\begin{equation}\label{eq}
R(T)_\alpha \equiv T + \{R(R(T)_\beta) \mid \beta \prec \alpha\}.
\end{equation}
More formally, let $\sigma_T(x)$ be an elementary numeration of a theory $T$.
Fix an elementary formula $\mathsf{Ax}_R(x, y)$ such that for each elementary formula $\tau(x)$ the formula $\mathsf{Ax}_R(x, \ul\tau\ur)$
numerates the schema $R(U)$, where $U$ is the theory numerated by $\tau$. W.l.o.g.\ we may also assume that
\begin{equation}\EA \vdash \forall x, y\, (\mathsf{Ax}_R(x, y) \imp x \geqslant y).\label{bd} \end{equation}
Consider an elementary formula $\rho(\alpha, x)$ defined by the following fixed-point equation:
\begin{equation}
\EA \vdash \forall \alpha\forall x\:(\rho(\alpha, x) \eqv (\sigma_T(x) \vee \exists \beta \prec \alpha\, \mathsf{Ax}_R(x, \ul \rho(\bar{\beta}, x) \ur))). \label{fpe}
\end{equation}
Note that $\beta \leqslant \ul \beta \ur \leqslant \ul \rho(\bar{\beta}, x) \ur$, whence the quantifier for $\beta$ can be bounded by $x$ in view of condition \refeq{bd}.
The formula $\rho(\alpha, x)$ thus constructed numerates a parametrized family of theories
that we denote $R(T)_\alpha$ and that provably satisfies equation (\ref{eq}).

We write $(T)^n_\alpha$ for $\RFN_{\Sigma_n}(T)_\alpha$ and $T_\alpha$ for $\Con(T)_\alpha$, the $\alpha$th member of
the progression of iterated consistency assertions.
The following lemma essentially due to Ulf Schmerl~\cite{Schm} provides a useful tool for reasoning about transfinite iterations inside weak theories such as $\EA$.
\begin{lemma}[reflexive induction]
For any elementary linear ordering $(D, \prec)$, any theory $T$ extending $\EA$ is closed under the following rule:
$$
\frac{\forall \alpha\, (\Box_T \forall \beta \prec \overline{\alpha}\, \varphi(\beta) \imp \varphi(\alpha))}{\forall \alpha\, \varphi(\alpha)}.
$$
\end{lemma}

Using this lemma it can be shown that the sequence of theories $(R(T)_\alpha)_{\alpha \in D}$ is unique modulo provable equivalence in $\EA$.
For more details on reflection principles and their transfinite iterations see~\cite{Bek99b}.

\section{Provable $n$-provability}
Given a natural number $n > 0$ and a pair of theories $T$ and $S$, define the following set of formulas
$$
C^n_S(T) := \{\varphi \mid T \vdash [n]_S\varphi\}.
$$
Since $[n]_S$ satisfies L{\"o}b's derivability conditions, $C^n_S(T)$ is a deductively closed set extending $S$.
Note that $C^n_S(T)$ can be viewed as a theory with the provability predicate $\Box_T[n]_S$.
In particular, $C^n_S(T) \subseteq U$ is a shorthand for $\forall \psi\, (\Box_T[n]_S\psi \imp \Box_U\psi)$.
Also notice that $T \subseteq U$ implies $C^n_S(T) \subseteq C^n_S(U)$.
We start with several simple facts about the theories $C^n_S(T)$.

The following basic result is due to Mingzhong Cai~\cite[Proposition 5.5]{Cai15} (unpublished), who considered the set of all formulas $\phi$ provably $n$-provable in a theory $T$ for some $n\in \nat$, that is, the theory $C^\infty_T(T):=\bigcup_{n \in \mathbb{N}} C^n_T(T)$. We include this result with a short direct proof.

\begin{proposition}[M.~Cai] \label{cai}
$
C^\infty_T(T) \equiv T + \mathsf{RFN}(T).
$
\end{proposition}
\begin{proof}
Since $T \vdash \varphi$ implies $T \vdash \Box_T \varphi$, we obtain $T \subseteq C^\infty_T(T)$.
Now, $\mathsf{RFN}(T)$ is equivalent to the theory axiomatized by $\{\langle n \rangle_T \top \mid n \in \mathbb{N}\}$.
For each $n$ we have
\begin{align*}
T \vdash [n]_T \bot &\rightarrow [n+1]_T\bot\\
&\rightarrow [n+1]_T\langle n \rangle_T \top.
\end{align*}
By provable $\Sigma_{n+2}$-completeness, since $\langle n \rangle_T\top$ is $\Pi_{n+1}$, we derive
$$
T \vdash \langle n \rangle_T \top \rightarrow [n+1]_T\langle n \rangle_T \top,
$$
whence $T \vdash [n+1]_T\langle n \rangle_T \top$ and $\langle n \rangle_T \top \in C^\infty_T(T)$.
This proves $T + \mathsf{RFN}(T) \subseteq C^\infty_T(T)$.

To show the converse, assume $T \vdash [n]_T \varphi$ for some $n \in \mathbb{N}$. Fix an arbitrary $m > n$ such that $\varphi \in \Pi_{m+1}$.
Using provable $\Sigma_{m+1}$-completeness we derive
\begin{align*}
T \vdash [n]_T \varphi \wedge \neg \varphi &\rightarrow [n]_T \varphi \wedge [m]_T \neg\varphi\\
&\rightarrow [m]_T \varphi \wedge [m]_T \neg\varphi\\
&\rightarrow [m]_T \bot,
\end{align*}
whence $T \vdash \langle m\rangle_T \top \rightarrow ([n]_T\varphi \rightarrow \varphi)$.
Combining it with $T \vdash [n]_T \varphi$, we obtain $T + \langle m \rangle_T \top \vdash \varphi$.
It follows that $T + \mathsf{RFN}(T) \vdash \varphi$, hence $C^\infty_T(T) \subseteq T + \mathsf{RFN}(T)$.
\end{proof}

In this paper we study a more delicate question of characterizing the theories $C^n_S(T)$ for a fixed $n>0$.
As a first simple observation we show that such theories are not sensitive to the extension of $T$ by axioms weaker than $\RFN_{\Pi_n}(S)$.

\begin{proposition}
For every sentence $\psi$ and every $\Pi_{n+1}$-sentence $\pi$,
\begin{itemize}
\item[(i)]  $C^n_S(T + \la n \ra_S \psi) \subseteq C^n_S(T) + \psi$.
\item[(ii)] $C^n_S(T + \pi) \subseteq C^n_S(T) + \pi$.
\end{itemize}
\end{proposition}
\begin{proof}
Note that (ii) follows from (i) by provable $\Sigma_{n+1}$-completeness. Indeed, since $\neg \pi$ is $\Sigma_{n+1}$,
we have $T \vdash \neg\pi \imp [n]_S\neg\pi$ or, equivalently,
$T + \la n\ra_S \pi \vdash \pi$, and by the monotonicity of $C^n_S(\cdot)$ operator and (i) it yields
$C^n_S(T + \pi) \subseteq C^n_S(T + \la n \ra_S \pi) \subseteq C^n_S(T) + \pi$.

Now let us prove (i). Assume $T + \la n \ra_S \psi \vdash [n]_S\varphi$ for some $\varphi$. Clearly, $T + \neg \la n \ra_S \psi \vdash [n]_S\neg\psi$.
Both $[n]_S\varphi$ and $[n]_S\neg\psi$ imply $[n]_S(\psi \imp \varphi)$, whence $T \vdash [n]_S(\psi \imp \varphi)$, that is,
$(\psi \imp \varphi) \in C^n_S(T)$. It follows that $C^n_S(T) + \psi \vdash \varphi$.
\end{proof}

\begin{corollary}\label{cor:rfn}
$C^n_S(T + \la n \ra_S \top)\equiv  C^n_S(T).$
\end{corollary}

It is well-known by the results of D.~Leivant~\cite{Lei83} and H.~Ono~\cite{Ono87} that
$\ISi_n \equiv \EA + \la n + 1\ra_\EA\top$ (cf also~\cite{Bek99b}). Hence we obtain
\begin{corollary}
$C^{n+1}_S(\EA) \equiv C^{n+1}_S(\ISi_n)$.
\end{corollary}

The following lemma shows that the theories $C^n_S(T)$ are strictly stronger than $S$.
\begin{lemma}\label{lm:rfn}
Provably in $\EA$, $S + \Rfn^n(S) \subseteq C^{n+1}_S(T)$.
\end{lemma}
\begin{proof}
Arguing in $\EA$ assume $S + \Rfn^n(S) \vdash \varphi$. We get a sequence of sentences $\psi_1, \dots, \psi_n$ such that
$$
\EA \vdash \Box_S \left(\bigwedge_{i = 1}^n ([n]_S \psi_i \imp \psi_i) \imp \varphi\right).
$$
Since $\EA \vdash \forall \psi\, (\Box_S\psi \imp [n+1]_S\psi)$, we have
$$
\EA \vdash [n+1]_S \left(\bigwedge_{i = 1}^n ([n]_S \psi_i \imp \psi_i) \imp \varphi\right).
$$
If we show that for each $i \in \{1, \dots, n\}$ it holds that $\EA \vdash [n+1]_S\left([n]_S \psi_i \imp \psi_i\right)$,
then from the previous derivation we obtain $\EA \vdash [n+1]_S\varphi$ and we get the required result since $T$ is an extension of $\EA$.
We derive
\begin{align*}
\EA \vdash [n]_S\psi_i \imp &\ [n+1]_S\psi_i \\
\imp &\ [n+1]_S([n]_S\psi_i \imp \psi_i)
\end{align*}
The sentence $\neg[n]_S\psi_i$ is $\Pi_{n+1}$, so we have
\begin{align*}
\EA \vdash \neg[n]_S\psi_i \imp &\ [n+1]_S\left(\neg[n]_S\psi_i\right) \\
\imp &\ [n+1]_S([n]_S\psi_i \imp \psi_i),
\end{align*}
whence $\EA \vdash [n+1]_S([n]_S\psi_i \imp \psi_i)$.
\end{proof}

It also answers the question concerning the axiomatization complexity of the theories $C^n_S(T)$.
By the Unboundedness theorem of Kreisel and L\'evy~\cite{KrL} (cf also~\cite[Corollary 2.22]{Bek05en}) we obtain
\begin{corollary}
For each $n > 0$ the theory $C^n_S(T)$ is of unbounded arithmetical complexity. In particular, it is not finitely axiomatized.
\end{corollary}

\section{Provable in $\EA$ $1$-provability}
Recall that the operator $C^n_S(\cdot)$ is monotone w.r.t.\ the inclusion of theories
and that we consider theories extending $\EA$.
So, the theory $C^1_S(\EA)$ is in some sense the base case for our investigations.
In this section we characterize this theory in terms of the local reflection principle over $S$.
For convenience we write $C_S(T)$ for $C^1_S(T)$ throughout the paper.
\begin{lemma}\label{lm:herbrand}
Provably in $\EA^+$, $C_S(\EA) \equiv S + \Rfn(S)$.
\end{lemma}
\begin{proof}

The inclusion $S + \Rfn(S) \subseteq C_S(\EA)$ follows from Lemma \ref{lm:rfn}.
We focus on the converse inclusion.
In this proof we will consider $\EA$ stated as a quantifier-free theory in the language with the terms for all elementary functions. Such a theory is known to be definitionally equivalent to the original formulation of $\EA$~(cf~\cite{Bek99b}).

Assume $\EA \vdash [1]_S\varphi$. Let $\forall x\, \delta(x, z)$, where $\delta(x, z)$ is a bounded formula,
be a $\Pi_1$-formula that is $\EA$-equivalent to $\True_{\Pi_1}(z)$.
By the definition of $[1]_S\varphi$ the assumption implies
$$
\EA \vdash \exists p\, \exists z\, \forall x\, (\delta(x, z) \wedge \Prf_S(y, \ul \forall x\, \delta(x, \overline{z}) \imp \varphi \ur)).
$$
We apply a version of Herbrand's theorem for $\Sigma_2$-formulas (which can be formalized in $\EA^+$)
to the derivation above. Thus, we obtain a sequence of terms
$t_0$, $p_0$, $t_1(x_1)$, $p_1(x_1), \dots, t_k(x_1, \dots, x_k)$, $p_k(x_1, \dots, x_k)$
such that the following disjunction is provable in $\EA$:
\begin{align*}
&\delta(x_1, t_0) \wedge \Prf_S(p_0, \ul \forall x\, \delta(x, \overline{t}_0) \imp \varphi \ur)\ \vee\\
&\delta(x_2, t_1(x_1)) \wedge \Prf_S(p_1(x_1), \ul \forall x\, \delta(x, \overline{t_1(x_1)}) \imp \varphi \ur)\ \vee\\
&\delta(x_3, t_2(x_1, x_2)) \wedge \Prf_S(p_2(x_1, x_2), \ul \forall x\, \delta(x, \overline{t_2(x_1, x_2)}) \imp \varphi \ur)\ \vee \\
&\dots \\
&\delta(x_{k+1}, t_k(x_1, \dots, x_k)) \wedge \Prf_S(p_k(x_1, \dots, x_k), \ul \forall x\, \delta(x, \overline{t_k(x_1, \dots, x_k)}) \imp \varphi \ur),
\end{align*}
with $x_1, x_2, \dots, x_{k+1}$ as free variables.
Our aim is to show that $\EA + \Rfn(S) \vdash \varphi$.
We do this by arguing informally in $\EA + \Rfn(S)$ and considering cases corresponding to the disjunction above.
Note that by provable $\Sigma_1$-completeness for each $\EA$-term $t(x_1, \dots, x_k)$ we have
$$
\EA \vdash \forall x_1\, \dots \forall x_k \left(\Box_S(\overline{t(x_1, \dots, x_k)} \eqv t(\overline{x}_1, \dots, \overline{x}_k))\right).
$$
This allows to replace occurrences of the form $\overline{t_i(x_1, \dots, x_i)}$ under $\Box_S$ with $t_i(\overline{x}_1, \dots, \overline{x}_i)$.

For $i \in \{1, \dots, k + 1\}$ let us denote by $C_i(x_1, \dots, x_i)$ the following elementary formula
$$
\delta(x_i, t_{i-1}(x_1, \dots, x_{i-1})) \wedge
\Prf_S(p_{i-1}(x_1, \dots, x_{i-1}), \ul \forall x\, \delta(x, \overline{t_{i-1}(x_1, \dots, x_{i-1})}) \imp \varphi \ur),
$$
that is, the $i$th member of the disjunction above.

We start with the first line by considering two cases: $\forall x_1\, C_1(x_1)$ and $\exists x_1\, \neg C_1(x_1)$.
Assume $\forall x_1\, C_1(x_1)$, then by the definition of $C_1(x_1)$ we have
$\forall x_1\, \delta(x_1, t_0)$ and $\Prf_S(p_0, \ul \forall x\, \delta(x, \overline{t}_0) \imp \varphi \ur)$.
The latter formula clearly implies $\Box_S(\forall x\, \delta(x, t_0) \imp \varphi)$.
Since $t_0$ is a closed term we can use $\Rfn(S)$ to get $\forall x\, \delta(x, t_0) \imp \varphi$,
hence we obtain $\varphi$.

Conversely, assume $\exists x_1 \neg C_1(x_1)$.
Using $\Delta_0(\exp)$-induction we can find the least element, denote it by $c_1$, satisfying $\neg C_1(x)$.
In other words, $c_1$ is the unique element satisfying the following elementary formula
$$
D_1(x) := \neg C_1(x) \wedge \forall y < x\, C_1(y),
$$
and clearly, $\EA \vdash \forall x, y\, (D_1(x) \wedge D_1(y) \imp x = y)$.
By provable $\Sigma_1$-completeness $D_1(c_1)$ implies $\Box_S D_1(\overline{c}_1)$.
Also by its definition $c_1$ falsifies the first line of the disjuction.

Now we substitute $c_1$ for $x_1$ in the disjunction above and consider two cases for the second line in the same way as for the first one.
Namely, assume $\forall x_2\, C_2(c_1, x_2)$, that is, $\forall x_2\, \delta(x_2, t_1(c_1))$ and
$\Prf_S(p_1(c_1), \ul \forall x\, \delta(x, \overline{t_1(c_1)}) \imp \varphi \ur)$.
The second formula implies $\Box_S(\forall x\, \delta(x, t_1(\overline{c}_1)) \imp \varphi)$.

We want to show
\begin{equation}\label{eq:h1}
\Box_S \left(\forall y\, (D_1(y) \imp (\forall x\, \delta(x, t_1(y)) \imp \varphi))\right).
\end{equation}
Indeed, arguing in $S$:
\begin{quote}
Assume $D_1(y)$. By $D_1(\overline{c}_1)$ and the uniqueness of $c_1$ we get $y = \overline{c}_1$.
But then $\forall x\, \delta(x, t_1(\overline{c}_1)) \imp \varphi$ implies $\forall x\, \delta(x, t_1(y)) \imp \varphi$, as required.
\end{quote}
Using $\Rfn(S)$ from \eqref{eq:h1} we obtain $\forall y\, (D_1(y) \imp (\forall x\, \delta(x, t_1(y)) \imp \varphi))$.
Instantiating $y = c_1$ we get $D_1(c_1) \imp (\forall x\, \delta(x, t_1(c_1)) \imp \varphi)$,
which, together with $D_1(c_1)$ and $\forall x\, \delta(x, t_1(c_1))$ implies $\varphi$.

Now, assume $\exists x_2 \neg C_2(c_1, x_2)$. We proceed in the same way as above.
Namely, by $\Delta_0(\exp)$-induction there exists the least $x$ satisfying $\neg C_2(c_1, x)$, denote it by $c_2$.
The pair $(c_1, c_2)$ satisfies the following elementary formula
$$
D_2(x, y) := \neg C_2(x, y) \wedge \forall z < y\, C_1(x, z),
$$
and as above we have $S \vdash \forall x, y, z\, (D_2(z, x) \wedge D_2(z, y) \imp x = y)$.
Together $c_1$ and $c_2$ falsify the first two lines of the disjunction.
Also by provable $\Sigma_1$-completeness $D_2(c_1, c_2)$ implies $\Box_S D_2(\overline{c}_1, \overline{c}_2)$.

We continue considering two cases for each line of the disjunction.
If $\forall x_i\, C_i(c_1, \dots, c_{i-1}, x_i)$, at first we eliminate the occurrences of
$c_1, \dots, c_{i-1}$ under $\Box_S$ by moving their elementary definitions under it as was done above.
Then we apply $\Rfn(S)$, obtaining a sentence of the form
$$
\forall x_1\, \dots \forall x_i\, (D_1(x_1) \wedge \dots \wedge D_i(x_1, \dots, x_i) \imp (\forall x\, \delta(x, t_i(x_1, \dots, x_i) \imp \varphi)),
$$
take $x_1 = c_1, \dots, x_i = c_i$ respectively, and derive $\varphi$ using hypotheses.

Conversely, if $\exists x_i \neg C_i(c_1, \dots, c_{i-1}, x_i)$, we define a new element $c_i$ as the
least $x$ satisfying $\neg C_i(c_1, \dots, c_{i-1}, x)$, which, together with previously defined $c_1, \dots, c_{i-1}$,
falsifies the first $i$ lines of the disjunction.

We proceed in this fashion by successively falsifying the lines. Since the whole disjunction is provable,
it cannot be that all the lines are falsified, hence for some $i \in \{1, \dots, k+1\}$
it must be the case that $\forall x_i\, C_i(c_1, \dots, c_{i-1}, x_i)$.
But in this case we derive $\varphi$, as required.
\end{proof}

\section{Provable $1$-provability for extensions of $\EA$}\label{sec:main}
In this section we study the case of an arbitrary theory $T$ and $n = 1$.
The main result of the section is Theorem \ref{th:main}, which characterizes the theory $C_S(T)$ in terms of a theory $S$ and the $\Sigma^0_2$-ordinal of a theory $T$.

We start by reducing the problem of axiomatizing $C_S(T)$ for an arbitrary arithmetical theory $T$
to the case of $\EA$, which has been dealt with in the previous section.
Note that since $[1]_S\varphi$ is a $\Sigma_2$-formula, the theory $C_S(T)$ depends only on the $\Sigma_2$-consequences of $T$.
$\Sigma_2$-consequences of the fragments of $\PA$ were studied in \cite{BekVis}.

Let us define $\omega_0 := \omega$, $\omega_{n+1} := \omega^{\omega_n}$ for $n < \omega$ and $\varepsilon_0 = \sup\{\omega_n \mid n < \omega\}$.
We fix some natural ordinal notation system for ordinals $\alpha < \varepsilon_0$ and its presentation in arithmetic.
Whenever we mention these ordinals in a formal context, we assume that this particular ordinal notation system is being used.

Assume an elementary well-ordering $(D, \prec)$ is fixed.
The \emph{$\Sigma^0_2$-ordinal} of a theory $T$ was defined in~\cite{BekVis} in terms of the iterations of the local $\Sigma_2$-reflection schema over $\EA$ as follows:
$$
|T|_{\Sigma^0_2} := \sup\{\alpha \in D \mid \Rfn_{\Sigma_2}(\EA)_\alpha \subseteq T\},
$$
where as usual we assume $|T|_{\Sigma^0_2} := \infty$ if the inclusion holds for all $\alpha \in D$.

A theory $T$ is said to be \emph{$\Sigma^0_2$-regular}, if there is an $\alpha \in D$ such that
$$
T \equiv_{\Sigma_2} \Rfn_{\Sigma_2}(\EA)_\alpha.
$$
Theories $\ISi_n$ were shown in \cite{BekVis} to be regular w.r.t.\ the natural ordinal notation system for $\varepsilon_0$.
\begin{lemma}\label{lm:sigma2-consequences}
$\Sigma_2$-consequences of $\ISi_n$ are axiomatized by $\Rfn_{\Sigma_2}(\EA)_{\omega_n}$ for $n > 0$.
\end{lemma}

It follows that $C_S(\ISi_n) \equiv C_S(\Rfn_{\Sigma_2}(\EA)_{\omega_n})$.
But also by \cite[Proposition 5.2]{Bek99b} we have $\forall \alpha\, \Rfn_{\Sigma_2}(T)_\alpha \equiv_{\Sigma_2} \Rfn(T)_\alpha$ provably in $\EA$,
whence $C_S(\ISi_n) \equiv C_S(\Rfn(\EA)_{\omega_n})$.
Therefore we now focus on characterizing theories of the form $C_S(\Rfn(\EA)_\alpha)$.

The following lemma shows that transfinite iterations of the local reflection principle along an arbitrary elementary well-ordering $(D, \prec)$
and the $C_S(\cdot)$ operator can be permuted.
This allows us to reduce the characterization of $C_S(\Rfn(\EA)_\alpha)$ to that of $C_S(\EA)$.
\begin{lemma}\label{lm:commutativity}
Provably in $\EA$, $\forall \alpha\, C_S(\Rfn(T)_\alpha) \equiv \Rfn(C_S(T))_\alpha$.
\end{lemma}
\begin{proof}
We give an informal argument by reflexive induction on $\alpha$ in $\EA$. The case $\alpha = 0$ is trivial, so we can assume $\alpha \neq 0$.
Denote $\Rfn(T)_\alpha$ and $\Rfn(C_S(T))_\alpha$ by $T^\alpha$ and $U^\alpha$ respectively.

We split the proof into two parts:
\begin{itemize}
\item[(i)] $\forall \alpha\, U^\alpha \subseteq C_S(T^\alpha).$
\item[(ii)] $\forall \alpha\, C_S(T^\alpha) \subseteq U^\alpha.$
\end{itemize}

(i) Assume $U^\alpha \vdash \varphi$.  By the definition of $U^\alpha$ and the formalized deduction theorem there exist sentences
$\psi_1, \dots, \psi_m$ and $\beta \prec \alpha$ such that
$$
U^\beta \vdash \bigwedge_{i = 1}^m (\Box_{U^\beta} \psi_i \rightarrow \psi_i) \rightarrow \varphi.
$$
By $\Sigma_1$-completeness we have
$$
\mathsf{EA} \vdash \Box_{U^\beta} \left(\bigwedge_{i = 1}^m (\Box_{U^\beta} \psi_i \rightarrow \psi_i) \rightarrow \varphi\right).
$$
The reflexive induction hypothesis for $\beta$ can be rewritten as
$$
\mathsf{EA} \vdash \forall \psi \left(\Box_{U^\beta} \psi \rightarrow \Box_{T^\beta}[1]_S\psi\right),
$$
hence
$$
\mathsf{EA} \vdash \Box_{T^\beta}[1]_S\left(\bigwedge_{i = 1}^m (\Box_{U^\beta} \psi_i \rightarrow \psi_i) \rightarrow \varphi\right).
$$
This implies
$$
T + \mathsf{Rfn}(T^\beta) \vdash [1]_S\left(\bigwedge_{i = 1}^m (\Box_{U^\beta} \psi_i \rightarrow \psi_i) \rightarrow \varphi\right).
$$
We will show that for each $i \in \{1, \dots, m\}$ we have
$$
T + \mathsf{Rfn}(T^\beta) \vdash [1]_S\left(\Box_{U^\beta} \psi_i \rightarrow \psi_i\right).
$$
Combining it with the previous derivation, we get $T^{\beta + 1} \vdash [1]_S\varphi$,
hence $T^\alpha \vdash [1]_S\varphi$, as required.

Fix some $i \in \{1, \dots, m\}$. Using the reflexive induction hypothesis we derive
\begin{align*}
T + \mathsf{Rfn}(T^\beta) \vdash \Box_{U^\beta}\psi_i \rightarrow &\ \Box_{T^\beta}[1]_S\psi_i \\
\rightarrow &\ \Box_{T^\beta}[1]_S(\Box_{U^\beta}\psi_i \rightarrow \psi_i) \\
\rightarrow &\ [1]_S(\Box_{U^\beta}\psi_i \rightarrow \psi_i).
\end{align*}
But since $\neg\Box_{U^\beta}\psi_i$ is a $\Pi_1$-sentence, we have
\begin{align*}
\EA \vdash \neg\Box_{U^\beta}\psi_i \rightarrow &\ [1]_S\left(\neg\Box_{U^\beta}\psi_i\right) \\
\rightarrow &\ [1]_S(\Box_{U^\beta}\psi_i \rightarrow \psi_i),
\end{align*}
that yields $T + \mathsf{Rfn}(T^\beta) \vdash [1]_S\varphi$.

(ii) Assume $T^\alpha \vdash [1]_S\varphi$. By the definition of $T^\alpha$ and the formalized deduction theorem there exist sentences
$\psi_1, \dots, \psi_n$ and $\beta \prec \alpha$ such that
\begin{equation} \label{eq:c1}
T^\beta \vdash \bigwedge_{i = 1}^m (\Box_{T^\beta} \psi_i \rightarrow \psi_i) \rightarrow [1]_S \varphi.
\end{equation}
Our aim is to show that for any subset $I \subseteq \{1, \dots, m\}$ it holds that
\begin{equation} \label{eq:c2}
U + \mathsf{Rfn}(U^\beta) \vdash \left(\bigwedge_{i \in I}\neg\Box_{T^\beta} \psi_i \wedge \bigwedge_{i \notin I} \Box_{T^\beta} \psi_i\right) \rightarrow \varphi,
\end{equation}
whence $U + \mathsf{Rfn}(U^\beta) \vdash \varphi$ follows by considering all $2^m$ cases, hence $U^\alpha \vdash \varphi$, as required.

Firstly, let us show that for any $I \subseteq \{1, \dots, m\}$ we have
\begin{equation} \label{eq:c3}
T^\beta \vdash \bigwedge_{i \in I} (\Box_{T^\beta} \psi_i \rightarrow \psi_i) \rightarrow [1]_S \varphi
\Longrightarrow T^\beta \vdash [1]_S \left(\left[\bigwedge_{i \in I}\neg\Box_{T^\beta} \psi_i\right] \rightarrow \varphi\right).
\end{equation}
Indeed, for each $i$ we have $T^\beta \vdash \neg\Box_{T^\beta} \psi_i \rightarrow (\Box_{T^\beta} \psi_i \rightarrow \psi_i)$,
so by the premise of \eqref{eq:c3} we get
\begin{align*}
T^\beta \vdash \left[\bigwedge_{i \in I}\neg\Box_{T^\beta} \psi_i\right] \rightarrow &
\bigwedge_{i \in I} (\Box_{T^\beta} \psi_i \rightarrow \psi_i) \\
\rightarrow &\ [1]_S \varphi \\
\rightarrow &\ [1]_S \left(\left[\bigwedge_{i \in I} \neg\Box_{T^\beta} \psi_i\right] \rightarrow \varphi\right).
\end{align*}
Using provable $\Sigma_1$-completeness for $[1]_S$ we derive
\begin{align*}
T^\beta \vdash \left[\bigvee_{i \in I} \Box_{T^\beta} \psi_i\right] \rightarrow &\
[1]_S \left(\bigvee_{i \in I} \Box_{T^\beta}\psi_i \right) \\
\rightarrow &\ [1]_S \left(\left[\bigwedge_{i \in I}\neg\Box_{T^\beta} \psi_i\right] \rightarrow \varphi\right),
\end{align*}
that yields the conclusion of \eqref{eq:c3}. Moreover, this argument can be formalized in $\mathsf{EA}$, so for each subset $I \subseteq \{1, \dots, m\}$ we have
\begin{equation} \label{eq:c4}
\mathsf{EA} \vdash \Box_{T^\beta} \left(\bigwedge_{i \in I} (\Box_{T^\beta} \psi_i \rightarrow \psi_i) \rightarrow [1]_S \varphi\right)
\rightarrow \Box_{T^\beta}[1]_S\left(\left[\bigwedge_{i \in I} \neg\Box_{T^\beta} \psi_i\right] \rightarrow \varphi\right).
\end{equation}

Now, fix some $I \subseteq \{1, \dots, m\}$. For each $i$ we have $T^\beta \vdash \psi_i \rightarrow (\Box_{T^\beta} \psi_i \rightarrow \psi_i)$, whence
$$
T^\beta \vdash \bigwedge_{i \notin I} \psi_i \rightarrow \bigwedge_{i \notin I}(\Box_{T^\beta} \psi_i \rightarrow \psi_i).
$$
Using this and \eqref{eq:c1} we derive
$$
T^\beta \vdash \bigwedge_{i \notin I} \psi_i \rightarrow \left(\bigwedge_{i \in I} (\Box_{T^\beta} \psi_i \rightarrow \psi_i)
\rightarrow [1]_S \varphi \right).
$$
By $\Sigma_1$-completeness this implies
$$
\mathsf{EA} \vdash \Box_{T^\beta} \bigwedge_{i \notin I} \psi_i
\rightarrow \Box_{T^\beta}\left(\bigwedge_{i \in I} (\Box_{T^\beta} \psi_i \rightarrow \psi_i)
\rightarrow [1]_S \varphi \right).
$$
We use \eqref{eq:c4} to obtain
$$
\mathsf{EA} \vdash \Box_{T^\beta} \bigwedge_{i \notin I} \psi_i
\rightarrow \Box_{T^\beta}[1]_S\left(\left[\bigwedge_{i \in I} \neg\Box_{T^\beta} \psi_i\right] \rightarrow \varphi\right).
$$
By the reflexive induction hypothesis for $\beta$ we have
$$
\mathsf{EA} \vdash \Box_{T^\beta}[1]_S\left(\left[\bigwedge_{i \in I} \neg\Box_{T^\beta} \psi_i\right] \rightarrow \varphi\right)
\rightarrow  \Box_{U^\beta}\left(\left[\bigwedge_{i \in I} \neg\Box_{T^\beta} \psi_i\right] \rightarrow \varphi\right),
$$
whence
$$
\mathsf{EA} \vdash \Box_{T^\beta} \bigwedge_{i \notin I} \psi_i
\rightarrow \Box_{U^\beta}\left(\left[\bigwedge_{i \in I} \neg\Box_{T^\beta} \psi_i\right] \rightarrow \varphi\right).
$$
Finally, this yields
$$
U + \mathsf{Rfn}(U^\beta) \vdash \Box_{T^\beta} \bigwedge_{i \notin I} \psi_i
\rightarrow \left(\left[\bigwedge_{i \in I} \neg\Box_{T^\beta} \psi_i\right] \rightarrow \varphi\right),
$$
and since
$$
\mathsf{EA} \vdash \left(\Box_{T^\beta} \bigwedge_{i \notin I} \psi_i\right) \leftrightarrow \left(\bigwedge_{i \notin I} \Box_{T^\beta} \psi_i\right)
$$
we obtain \eqref{eq:c2}.
\end{proof}



Now we are ready to state and prove a general result characterizing the theory $C_S(T)$
in terms of the $\Sigma^0_2$-ordinal of $T$ and iterated local reflection over $S$.
\begin{theorem}\label{th:main}
If $T$ is a $\Sigma^0_2$-regular theory with $|T|_{\Sigma^0_2} = \alpha$, then $C_S(T) \equiv \Rfn(S)_{1 + \alpha}$.
\end{theorem}
\begin{proof}
By the discussion preceding Lemma \ref{lm:commutativity} and the hypothesis we have
$$
C_S(T) \equiv C_S(\Rfn_{\Sigma_2}(\EA)_\alpha) \equiv C_S(\Rfn(\EA)_\alpha).
$$
Using Lemmas \ref{lm:herbrand} and \ref{lm:commutativity} we obtain the following chain of equivalences
$$
C_S(T) \equiv C_S(\Rfn(\EA)_\alpha) \equiv \Rfn(C_S(\EA))_\alpha \equiv \Rfn(S + \Rfn(S))_\alpha \equiv \Rfn(S)_{1 + \alpha}.
$$
\end{proof}

Note that the equivalences stated in Lemmas \ref{lm:herbrand} and \ref{lm:commutativity} are provable in $\EA^+$.
It follows that if the theory is $\Sigma^0_2$-regular provably in $\EA^+$, then the conclusion of Theorem \ref{th:main} is also provable in $\EA^+$.

Now we obtain several corollaries characterizing the theories $C_S(T)$ for various fragments of $\PA$.

\begin{corollary}\label{cor:isigma}
For all $n > 0$ provably in $\EA^+$ we have $C_S(\ISi_n) \equiv \Rfn(S)_{\omega_n}$.
\end{corollary}
\begin{proof}
Due to the Lemma \ref{lm:sigma2-consequences} we can apply Theorem \ref{th:main} with $\alpha = \omega_n$.
We get the required result since $1 + \omega_n = \omega_n$ for $n > 0$.
It can be seen that the proof of Lemma \ref{lm:sigma2-consequences} can be formalized in $\EA^+$, whence
the result about $\EA^+$-provability of the equivalence follows.
\end{proof}

In particular, we answer the question concerning $C^1_\PA(\PA)$ raised in the introduction.
\begin{corollary}
Provably in $\EA^+$, $C_S(\PA) \equiv \Rfn(S)_{\varepsilon_0}$.
\end{corollary}

Using well-known conservation results we can also characterize the theories $C_S(T)$ for parameter free induction schemata.
\begin{corollary}\label{cor:parfree}
For all $n > 0$ we have $C_S(\ISi^-_n) \equiv C_S(\IPi^-_{n+1}) \equiv \Rfn(S)_{\omega_n}$.
\end{corollary}
\begin{proof}
By \cite[Proposition 7.6]{Bek99b} we have the following conservation results
$$
\ISi_n \equiv_{\Sigma_{n+2}} \ISi^-_n \equiv_{\mathcal{B}(\Sigma_{n+1})} \IPi^-_{n+1}.
$$
Since $n > 0$, these three theories have the same $\Sigma_2$-consequences,
hence
$$
C_S(\ISi_n) \equiv C_S(\ISi^-_n) \equiv C_S(\IPi^-_{n+1}).
$$
The result now follows from Corollary \ref{cor:isigma}.
\end{proof}

The next corollary covers the exceptional case of the theory $\IPi^-_1$.
\begin{corollary}\label{cor:ipi1}
$C_S(\IPi_1^-) \equiv \Rfn(S)_2$.
\end{corollary}
\begin{proof}
By \cite[Theorem 3]{Bek99} we have $\EA^+ + \IPi_1^- \equiv \EA^+ + \Rfn_{\Sigma_2}(\EA)$.
Since  $\EA^+ \equiv \EA + \la 1 \ra_\EA \top$, Corollary \ref{cor:rfn} implies that
$$
C_S(\IPi_1^-) \equiv C_S(\EA^+ + \IPi_1^-) \equiv C_S(\EA^+ + \Rfn_{\Sigma_2}(\EA)) \equiv C_S(\EA + \Rfn_{\Sigma_2}(\EA)) \equiv \Rfn(S)_2,
$$
where the last equivalence is due to Theorem \ref{th:main}.
\end{proof}

The following corollary shows how to compute the $\Pi^0_1$-ordinal for the theories $C_\EA(T)$ from the $\Sigma^0_2$-ordinal of $T$.
\begin{corollary}
If $T$ is a $\Sigma^0_2$-regular theory $|T|_{\Sigma^0_2} = \alpha$, then $|C_\EA(T)|_{\Pi^0_1} = \omega^{1 + \alpha}$.
\end{corollary}
\begin{proof}
By \cite[Proposition 6.2]{Bek99b} we have $\forall \beta\, \Rfn(\EA)_\beta \equiv_{\Pi_1} \EA_{\omega^\beta}$.
Applying Theorem \ref{th:main} and this fact with $\beta = 1 + \alpha$ we get
$$
C_\EA(T) \equiv \Rfn(\EA)_{1+\alpha} \equiv_{\Pi_1} \EA_{\omega^{1+\alpha}}.
$$
\end{proof}

\section{Relativization}
In this section we generalize the results obtained in previous sections to the notion of $n$-provability for $n > 1$.
The results concerning $\Sigma_2$-conservativity mentioned in the discussion before Lemma \ref{lm:commutativity}
have the corresponding relativizations to the case of $\Sigma_{n+2}$-formulas.
In particular, by the results of \cite{BekVis} we have the following
\begin{lemma}\label{lm:rel-sigma2-consequences}
$\Sigma_{n+2}$-consequences of $\ISi_m$ are axiomatized by $\Rfn^n_{\Sigma_{n+2}}(\EA)_{\omega_{m - n}}$ for $0 \leqslant n < m$.
\end{lemma}

By the relativization of \cite[Proposition 5.2]{Bek99b} we have $\forall \alpha\, \Rfn^n_{\Sigma_{n+2}}(T)_\alpha \equiv_{\Sigma_{n+2}} \Rfn^n(T)_\alpha$ provably in $\EA$.
Since $[n+1]S_\varphi$ is a $\Sigma_{n+2}$-formula, we obtain the analogous chain of equalities
$$
C^{n+1}_S(\ISi_m) \equiv C^{n+1}_S(\Rfn^n_{\Sigma_{n+2}}(\EA)_{\omega_{m - n}}) \equiv C^{n+1}_S(\Rfn^n(\EA)_{\omega_{m - n}}).
$$
These facts allow us to use the same line of the argument as in the previous sections.
We focus on the theories of the form $C^{n+1}_S(\Rfn^n(\EA))$ and reduce their characterization to that of $C^{n+1}_S(\EA)$.

The following lemma is the relativized version of Lemma \ref{lm:commutativity}.
\begin{lemma}\label{lm:rel-commutativity}
Provably in $\EA$, $C^{n+1}_S(\Rfn^n(T)_\alpha) \equiv \Rfn^n(C^{n+1}_S(T))_\alpha.$
\end{lemma}
\begin{proof}
The proof is obtained by the straightforward relativization of the proof of Lemma \ref{lm:commutativity} but with minor modifications.
Namely, one need to show that the reflexive induction hypothesis implies its relativized versions:
\begin{align}
\label{eq:rc1}
\EA &\vdash  \forall \psi \left([n]_{U^\beta} \psi \imp [n]_{T^\beta}[n+1]_S\psi\right),\\
\label{eq:rc2}
\EA &\vdash  \forall \psi \left([n]_{T^\beta}[n+1]_S\psi \imp [n]_{U^\beta}\psi\right).
\end{align}

We argue informally in $\EA$. To prove \eqref{eq:rc1} assume $U^\beta \vdash \pi \imp \psi$ for some true $\Pi_n$-sentence $\pi$.
Then by the reflexive induction hypothesis $T^\beta \vdash [n+1]_S(\pi \imp \psi)$,
and hence $T^\beta \vdash [n+1]_S\pi \imp [n+1]_S\psi$. But by provable $\Sigma_{n+2}$-completeness
we have $\EA \vdash \pi \imp [n+1]_S\pi$, which implies $T^\beta \vdash \pi \imp [n+1]_S\psi$, as required.

To prove \eqref{eq:rc2} assume $T^\beta \vdash \pi \imp [n+1]_S \psi$ for some true $\Pi_n$-sentence $\pi$.
By $\Sigma_{n+2}$-completeness we have $\EA \vdash \neg\pi \imp [n+1]_S\neg\pi$. Both $\psi$ and $\neg \pi$
imply $\pi \imp \psi$ and hence $T^\beta \vdash [n+1]_S(\pi \imp \psi)$. By the reflexive induction hypothesis
we then obtain $U^\beta \vdash \pi \imp \psi$.
\end{proof}

We present two proofs of the relativized version of Lemma \ref{lm:herbrand}.
Firstly, we give a model-theoretic proof.
\begin{lemma}\label{lm:mod-th}
$C^{n+1}_S(\EA) \subseteq \EA + \Rfn^n(S)$.
\end{lemma}
\begin{proof}
Assume $\EA \vdash [n+1]_S\varphi$ for some $\varphi$.
We will show that $\varphi$ is true in every model of $\EA + \Rfn^n(S)$.
Fix an arbitrary $M \models \EA + \Rfn^n(S)$ and note that $M \models \ISi^-_n$ since $\ISi^-_n \equiv \EA + \Rfn^n_{\Sigma_{n+1}}(\EA)$ by \cite[Theorem 1]{Bek99b}.
It is known that if $M \models \ISi^-_n$, then $K^{n+1}(M) \prec_{\Sigma_{n+1}} M$ for $n \geqslant 0$ (see Remark (i) after \cite[Theorem 2.1]{KPD}),
and also $K^{n+1}(M) \models \EA$, where $K^{n+1}(M)$ is the substructure of $M$ consisting of all $\Sigma_{n+1}$-definable elements (without parameters).

Since $\EA \vdash [n+1]_S\varphi$ it follows that $K^{n+1}(M) \models [n+1]_S\varphi$, that is, $\exists a \in K^{n+1}(M)$ such that
$$
K^{n+1}(M) \models \True_{\Pi_{n+1}}(a) \wedge \Box_S(\True_{\Pi_{n+1}}(\overline{a}) \imp \varphi).
$$
The formula above is $\Pi_{n+1} \wedge \Sigma_1$ and hence by $K^{n+1}(M) \prec_{\Sigma_{n+1}} M$ we have
\begin{equation}\label{eq:mtp1}
M \models \True_{\Pi_{n+1}}(a) \wedge \Box_S(\True_{\Pi_{n+1}}(\overline{a}) \imp \varphi),
\end{equation}
Fix some $\Sigma_{n+1}$-formula $\sigma(x)$ defining $a$, that is, $M \models \exists ! x \, \sigma(x) \wedge \sigma(a)$.
Assume $\sigma(x) = \exists y\, \pi(x, y)$, where $\pi(x, y)$ is $\Pi_n$, and consider the following formula
$$
\delta(z) := \pi((z)_0, (z)_1) \wedge \forall y < z \, \neg\pi((y)_0, (y)_1),
$$
where $(z)_0$ and $(z)_1$ are the components of the pair coded by $z$. By $M \models \exists x\, \sigma(x)$ we obtain $M \models \exists z\, \pi((z)_0, (z)_1)$.
But since $\pi((z)_0, (z)_1)$ is $\Pi_n$ we get $M \models \exists z\, \delta(z)$ by using $\mathsf{L}\Pi^-_n$ (due to \cite[Proposition 1.4]{KPD}).

We fix some $c \in M$ such that $M \models \delta(c)$ and note that, in particular, $M \models \pi((c)_0, (c)_1)$,
and hence $M \models \sigma((c)_0)$ by the definition of $\sigma(x)$.
But then $M \models (c)_0 = a$ since $\sigma(x)$ defines $a$.
Thus, \eqref{eq:mtp1} can be rewritten as
\begin{equation}\label{eq:mtp2}
M \models \True_{\Pi_{n+1}}((c)_0) \wedge \Box_S(\True_{\Pi_{n+1}}((\overline{c})_0) \imp \varphi).
\end{equation}
In case $n = 0$ the formula $\delta(z)$ is bounded and we take $\gamma(z) := \delta(z)$.
If $n > 0$ denote by $\tilde{\delta}(z)$ the following formula
$$
\pi((z)_0, (z)_1) \wedge \exists u\, \forall y < z \, \exists v < u\, \psi(y, v),
$$
where $\neg\pi((y)_0, (y)_1)$ is $\exists u\, \psi(y, u)$ with $\psi(u, y)$ in $\Pi_{n-1}$.
This formula is seen to be $\EA$-equivalent to some $\Sigma_{n+1}$-formula, that is,
there exists $\gamma(z)$ in $\Sigma_{n+1}$ such that $\EA \vdash \forall z\, (\gamma(z) \eqv \tilde{\delta}(z))$.
Using the corresponding collection axiom as in \cite[Proposition 1.7]{KPD} we obtain $M \models \forall z\, (\tilde{\delta}(z) \eqv \delta(z))$
and hence $M \models \forall z\, (\gamma(z) \eqv \delta(z))$.
In particular, $M \models \gamma(c)$.

Since $\gamma(z)$ is in $\Sigma_{n+1}$ we get $M \models [n]_S \gamma(\overline{c})$ by provable $\Sigma_{n+1}$-completeness
and hence $M \models [n]_S \delta(\overline{c})$, because by logic we have $\EA \vdash \forall z\, (\tilde{\delta}(z) \imp \delta(z))$
and so $\EA \vdash \forall z\, (\gamma(z) \imp \delta(z))$. Also by the definition of $\delta(z)$ we clearly have
\begin{equation}\label{eq:mtp3}
\EA \vdash \forall x, y\, (\delta(x) \wedge \delta(y) \imp x = y).
\end{equation}

Let us show that
\begin{equation}\label{eq:mtp4}
M \models [n]_S \left(\forall y\, (\delta(y) \imp (\True_{\Pi_{n+1}}((y)_0) \imp \varphi))\right).
\end{equation}
Indeed, arguing under $[n]_S$ inside $M$ assume $\delta(y)$. Since we also have $\delta(\overline{c})$, using \eqref{eq:mtp3}
we conclude that $y = \overline{c}$. But by \eqref{eq:mtp2} we have $\True_{\Pi_{n+1}}((\overline{c})_0) \imp \varphi$ and hence
$\True_{\Pi_{n+1}}((y)_0) \imp \varphi$, as required. We apply $\Rfn^n(S)$ to \eqref{eq:mtp4} and obtain
$$
M \models \forall y\, (\delta(y) \imp (\True_{\Pi_{n+1}}((y)_0) \imp \varphi)).
$$
In particular, we have $M \models \delta(c) \imp (\True_{\Pi_{n+1}}((c)_0) \imp \varphi)).$
Combining it with \eqref{eq:mtp2} and $M \models \delta(c)$ we get $M \models \varphi$.
\end{proof}

In order to have the result provable in $\EA^+$ as in the case of $1$-provability
we also give a syntactic proof that can be formalized in $\EA^+$.

\begin{lemma}\label{lm:rel-herbrand}
Provably in $\EA^+$, $C^{n+1}_S(\EA) \subseteq \ISi_n + \Rfn^n(S)$.
\end{lemma}
\begin{proof}
The proof is similar to that of Lemma \ref{lm:herbrand} and is based on the application of Herbrand's theorem.
In order to apply Herbrand's theorem we reduce the problem to the case of $\Sigma_2$-formulas by adapting the method from \cite[Theorem 4]{Bek97a}.
At first, let us consider the case $n = 1$. Assume $\varphi \in C^2_S(\EA)$, that is, $\EA \vdash [2]_S\varphi$.
Let $\forall x\, \exists y\, \delta(x, y, z)$, where $\delta(x, y, z)$ is a bounded formula,
be a $\Pi_2$-formula that is $\EA$-equivalent to $\True_{\Pi_2}(z)$, then by the definition of $[2]_S$ we have
\begin{equation}\label{eq:rh1}
\EA \vdash \exists z\, (\forall x\, \exists y\, \delta(x, y, z) \wedge \Box_S(\forall x\, \exists y\, \delta(x, y, \overline{z}) \imp \varphi)).
\end{equation}
Our aim is to eliminate the innermost existential quantifier by introducing a new function symbol to the language of $\EA$.
Let us consider the theory $U$ obtained from $\EA^g$, where $g$ is a new function symbol, by adding the following axiom
\begin{equation}\label{eq:rh2}
\forall x\, \forall z \leqslant x\, (\exists y\, \delta((z)_0, y, (z)_1) \imp \exists y \leqslant g(x) \, \delta((z)_0, y, (z)_1)).
\end{equation}

By \cite[Lemma 9.4]{Bek97a} there exists a non-relativizing interpretation $(\cdot)^-$ of $U$ in $\ISi_1$
that is identical on formulas in the language of $\EA$.
It follows that (provably in $\EA$) $U$ is conservative over $\ISi_1$.

The existential quantifier in the formula $\exists y\, \delta(x, y, z)$ can be bounded in $U$ using axiom \eqref{eq:rh2},
and the function that gives the least witness for this quantifier is elementary in $g$.
It follows that this function can be defined by some $\EA^g$-term $g'(x, z)$ such that
\begin{equation}\label{eq:rh3}
U \vdash \exists y\, \delta(x, y, z) \eqv \delta(x, g'(x, z), z).
\end{equation}
Using \eqref{eq:rh1} and \eqref{eq:rh3} we get
$$
U \vdash \exists p\, \exists z\, \forall x\, (\delta(x, g'(x, z), z) \wedge \Prf_S(p, \ul \forall x\, \exists y\, \delta(x, y, \overline{z}) \imp \varphi \ur)).
$$
Now, since $U$ has $\Pi^g_1$-axiomatization and the formula above is $\Sigma^g_2$ we are able to apply a version of Herbrand's theorem for $\Sigma_2$-formulas.
Let us denote by $C(x, z, p)$ the following formula
$$
\delta(x, g'(x, z), z) \wedge \Prf_S(p, \ul \forall x\, \exists y\, \delta(x, y, \overline{z}) \imp \varphi \ur))).
$$
By Herbrand's theorem we get a sequence of $\EA^g$-terms
$p_0$, $t_0$, $p_1(x_1)$, $t_1(x_1)$, \dots, $p_k(x_1, \dots, x_k)$, $t_k(x_1, \dots, x_k)$ such that the following disjunction is provable in $U$:
$$
C(x_1, t_0, p_0) \vee C(x_2, t_1(x_1), p_1(x_1)) \vee \dots \vee C(x_{k+1}, t_k(x_1, \dots, x_k), p_k(x_1, \dots, x_k))
$$
with $x_1, x_2, \dots, x_{k+1}$ as free variables.

By applying \eqref{eq:rh3}, weakening $\Prf_S$ to $\Box_S$ and rewriting occurrences of terms $t_0, t_1(x_1)$, $\dots$, $t_k(x_1, \dots, x_k)$
using existential quantifiers we obtain
\begin{align*}
U \vdash &\exists u_0\, \forall x_1\, \Bigl(\left(\exists y\, \delta(x_1, y, u_0) \wedge u_0 = t_0 \wedge
\Box_S (\forall x\, \exists y\, \delta(x, y, \overline{u}_0) \imp \varphi)\right)\ \vee \\
&\quad \exists u_1 \forall x_2 \Bigl( \left(\exists y\, \delta(x_2, y, u_1) \wedge u_1 = t_1(x_1) \wedge
\Box_S (\forall x\, \exists y\, \delta(x, y, \overline{u}_1) \imp \varphi)\right)\ \vee \\
&\qquad \qquad \dots \\
&\qquad \exists u_k \forall x_{k+1} \Bigl( \left(\exists y\, \delta(x_{k+1}, y, u_k) \wedge u_k = t_k(x_1, \dots, x_k) \wedge
\Box_S (\forall x\, \exists y\, \delta(x, y, \overline{u}_k) \imp \varphi)\right)\Bigr) \dots \Bigr).
\end{align*}
Let us consider the case $k = 1$. We apply the interpretation and obtain
\begin{align*}
\ISi_1 \vdash \exists u_0\, \forall x_1\, \Bigl(& \left(\exists y\, \delta(x_1, y, u_0) \wedge (u_0 = t_0)^- \wedge
\Box_S (\forall x\, \exists y\, \delta(x, y, \overline{u}_0) \imp \varphi)\right)\ \vee \\
&\exists u_1 \forall x_2 \Bigl( \left(\exists y\, \delta(x_2, y, u_1) \wedge (u_1 = t_1(x_1))^- \wedge
\Box_S (\forall x\, \exists y\, \delta(x, y, \overline{u}_1) \imp \varphi)\right)\Bigr)\Bigr).
\end{align*}
Note that all $(u_i = t_i(x_1, \dots, x_i))^-$ are equivalent to $\Sigma_2$-formulas in $\ISi_1$.
By induction on the complexity of a $\EA^g$-term $t(x_1, \dots, x_m)$ using \cite[Lemma 9.4]{Bek97a} for the base case
one can show
\begin{equation}\label{eq:rh4}
\ISi_1 \vdash \forall x_1\, \dots \forall x_m\, \exists ! u\, (u = t(x_1, \dots, x_m))^-.
\end{equation}
In particular, this holds for the formulas $(u_i = t_i(x_1, \dots, x_i))^-$.

Now we argue inside $\ISi_1 + \Rfn^1(\ISi_1 + S)$ as in the proof of Lemma \ref{lm:herbrand}
by considering cases corresponding to the disjunction above.

Fix some witness $u_0$ to the first existential quantifier and start from the first line.
In case $\neg ((u_0 = t_0)^- \wedge \Box_S (\forall x\, \exists y\, \delta(x, y, \overline{u}_0) \imp \varphi))$
the first line is falsified by any $x_1$, so we fix $x_1 = \bar{0}$ and move to the second line.
Note that we are happen to be in the same position as for the first line, and we can apply the same reasoning.

So assume $(u_0 = t_0)^- \wedge \Box_S (\forall x\, \exists y\, \delta(x, y, \overline{u}_0) \imp \varphi)$.
In this case we also consider two possibilities. Suppose $\forall x_1\, \exists y\, \delta(x_1, y, u_0)$.
By provable $\Sigma_2$-completeness we get $[1]_{\ISi_1 + S}(\overline{u}_0 = t_0)^-$.
Also by the assumption we have $\Box_S (\forall x\, \exists y\, \delta(x, y, \overline{u}_0) \imp \varphi)$.

We claim that $[1]_{\ISi_1 + S}\left(\forall v_0\, \left((v_0 = t_0)^- \imp (\forall x\, \exists y\, \delta(x, y, v_0) \imp \varphi)\right)\right)$.
We argue under $[1]_{\ISi_1 + S}$:
\begin{quote}
Assume $(v_0 = t_0)^-$. By $(\overline{u}_0 = t_0)^-$ and \eqref{eq:rh4} we get $v_0 = \overline{u}_0$,
whence \mbox{$\forall x\, \exists y\, \delta(x, y, v_0) \imp \varphi$} since $\forall x\, \exists y\, \delta(x, y, \overline{u}_0) \imp \varphi$.
\end{quote}
Applying $\Rfn^1(\ISi_1 + S)$ we get $\forall v_0\, \left((v_0 = t_0)^- \imp (\forall x\, \exists y\, \delta(x, y, v_0) \imp \varphi)\right)$.
By taking $v_0 = u_0$ and using hypotheses we obtain $\varphi$.

Conversely, assume $\exists x_1\, \forall y\, \neg \delta(x_1, y, u_0)$.
Let $B(u_0, x_1)$ be the formula asserting that $x_1$ is the least element satisfying $\forall y\, \neg \delta(x_1, y, u_0)$.
Using $\mathsf{L}\Pi_1$ we find such $x_1$ and denote it by $c_1$. Note that $B(u_0, x_1)$ is of the form $\Pi_1 \wedge \forall^b \Sigma_1$,
hence it is $\ISi_1$-equivalent to $\Sigma_2$-formula by using $\BSi_1$. Also we clearly have
$$
\ISi_1 \vdash \forall u\, \forall x\, \forall y\, (B(u, x) \wedge B(u, y) \imp x = y).
$$
The elements $u_0$ and $c_1$ together falsify the first line. We fix them and move to the second line.

Fix a witness $u_1$ for the second existential quantifier (with $x_1$ replaced by $c_1$).
In our case ($k = 1$) it must be that
\begin{equation}\label{eq:rh5}
\forall x_2\, \exists y\, \delta(x_2, y, u_1) \wedge (u_1 = t_1(c_1))^- \wedge \Box_S (\forall x\, \exists y\, \delta(x, y, \overline{u}_1) \imp \varphi).
\end{equation}
Using provable $\Sigma_2$-completeness we obtain
$[1]_{\ISi_1 + S}((\overline{u}_0 = t_0)^- \wedge B(\overline{u}_0, \overline{c}_1) \wedge (\overline{u}_1 = t_1(\overline{c}_1))^-).$

We claim that
$$
[1]_{\ISi_1 + S}\left(\forall v_0\, \forall v_1\, \forall d_1 \left( (v_0 = t_0)^-  \wedge B(v_0, d_1) \wedge (v_1 = t_1(d_1))^- \imp
(\forall x\, \exists y\, \delta(x, y, v_1) \imp \varphi)\right)\right).
$$
We argue under $[1]_{\ISi_1 + S}$ as follows:
\begin{quote}
Assume $(v_0 = t_0)^-  \wedge B(v_0, d_1) \wedge (v_1 = t_1(d_1))^-$.
The first hypothesis implies $v_0 = \overline{u}_0$, whence the second hypothesis implies
$d_1 = \overline{c}_1$, so $v_1 = \overline{u}_1$ by the third hypothesis.
By the assumption \eqref{eq:rh5} we get $\forall x\, \exists y\, \delta(x, y, v_1) \imp \varphi$.
\end{quote}
Using $\Rfn^1(\ISi_1 + S)$, taking $v_0 = u_0$, $v_1 = u_1$, $d_1 = c_1$ and applying the hypotheses we derive $\varphi$.

The proof for $k > 1$ is no different. It goes along the same lines as the proof of Lemma \ref{lm:herbrand}
but with the appropriate modifications as it was done above for the case $k = 1$.

Finally, since $\ISi_1$ is finitely axiomatizable we have $\ISi_1 + \Rfn^1(S) \vdash \Rfn^1(\ISi_1 + S)$.
Indeed, denote by $\psi$ the conjunction of all axioms of $\ISi_1$. By the formalized deduction theorem we get
\begin{align*}
\ISi_1 + \Rfn^1(S) \vdash [1]_{\ISi_1 + S}\theta &\imp [1]_S(\psi \imp \theta)\\
&\imp (\psi \imp \theta),
\end{align*}
and we obtain $\ISi_1 + \Rfn^1(S) \vdash [1]_{\ISi_1 + S}\theta \imp \theta$ since $\ISi_1 \vdash \psi$.
Hence $\ISi_1 + \Rfn^1(S) \vdash \varphi$.

A proof for an arbitrary $n > 1$ is essentially the same, but now we use \cite[Lemma 9.6]{Bek97a} and
\cite[Lemma 9.7]{Bek97a} to get the interpretation $(\cdot)^-$ and to show that $\True_{\Pi_{n+1}}(z)$ is equivalent to
some $\Pi^{g_1, \dots, g_n}_1$-formula. The formulas $(u_i = t_i(x_1, \dots, x_i))^-$ are equivalent to $\Sigma_{n+1}$-formulas in $\ISi_n$
and hence we get $[n]_{\ISi_n + S}$ instead of $[1]_{\ISi_1 + S}$ in the last part of the argument because
we apply provable $\Sigma_{n+1}$-completeness instead of $\Sigma_2$-completeness.
\end{proof}

Note that, while Lemmas \ref{lm:rfn} and \ref{lm:mod-th} give the equivalence $C^{n+1}_S(\EA) \equiv S + \Rfn^n(S)$
for all $S$ extending $\EA$, in order to have it provable in $\EA^+$ we need $S$ to be as strong as $\ISi_n$ due to the Lemma \ref{lm:rel-herbrand}.
Hence we get the following analog of Lemma \ref{lm:herbrand}.
\begin{lemma}\label{lm:rel-ea}
$C^{n+1}_S(\EA) \equiv S + \Rfn^n(S)$. Moreover, if $S$ extends $\ISi_n$ this equivalence is provable in $\EA^+$.
\end{lemma}

Using Lemmas \ref{lm:rel-commutativity} and \ref{lm:rel-ea} we obtain the relativizations of the main results of the previous section.
The notion of the $\Sigma^0_{n+2}$-ordinal of a theory $T$ and of a $\Sigma^0_{n+2}$-regular theory are defined analogously.
\begin{theorem}\label{th:rel-main}
If $T$ is a $\Sigma^0_{n+2}$-regular theory with $|T|_{\Sigma^0_{n+2}} = \alpha$, then $C^{n+1}_S(T) \equiv \Rfn^n(S)_{1 + \alpha}$.
\end{theorem}

Again if $\Sigma^0_{n+2}$-regularity of $T$ is provable in $\EA^+$ and $S$ extends $\ISi_n$,
the equivalence stated in Theorem \ref{th:rel-main} is provable in $\EA^+$.
This theorem together with Lemma \ref{lm:rel-sigma2-consequences} imply the following
\begin{corollary}
For all $0 \leqslant n < m$ we have $C^{n+1}_S(\ISi_m) \equiv \Rfn(S)_{\omega_{m - n}}$.
\end{corollary}

\begin{corollary}
$C^{n+1}_S(\PA) \equiv \Rfn^n(S)_{\varepsilon_0}$.
\end{corollary}

We also obtain the analogues of Corollaries \ref{cor:parfree} and \ref{cor:ipi1}.
\begin{corollary}
For all $0 \leqslant n < m$ we have $C^{n+1}_S(\ISi^-_m) \equiv C^{n+1}_S(\IPi^-_{m+1}) \equiv \Rfn(S)_{\omega_{m-n}}$.
\end{corollary}

\begin{corollary}
$C^{n+1}_S(\IPi_{n+1}^-) \equiv \Rfn^n(S)_2$ for $n > 0$.
\end{corollary}
\begin{proof}
By \cite[Theorem 1]{Bek99} we have $\IPi_{n+1}^- \equiv \EA + \Rfn^n_{\Sigma_{n+2}}(\EA)$ for $n > 0$. The result now follows by Theorem \ref{th:rel-main}.
\end{proof}

\section{Speed-up results}
Theorem \ref{th:main} implies that the theories $C_\EA(\EA)$ and $\EA + \Rfn(\EA)$ are deductively equivalent and this equivalence is provable in $\EA^+$.
In this section we prove a speed-up result, which shows that it is in some sense easier to prove
$1$-provability of a formula in $\EA$ rather than to find a proof of this formula in $\EA + \Rfn(\EA)$.

Namely, we show that $C_\EA(\EA)$ has \emph{superexponential speed-up} over $\EA + \Rfn(\EA)$, that is,
there is a sequence of sentences $\varphi_n$, a polynomial $p(n)$ and a constant $\varepsilon > 0$ such that for each $n$ there is
an $\EA$-proof of $[1]_\EA\varphi_n$ of size $p(n)$ and there is no $\EA+\Rfn(\EA)$-proof of $\varphi_n$ of size less than $(2^1_n)^\varepsilon$.
In particular, there is no elementary function $f(x)$ such that $f(p)$ codes an $\EA + \Rfn(\EA)$-proof of $\varphi$,
whenever $p$ codes an $\EA$-proof of $[1]_\EA \varphi$ for some $\varphi$.


Our strategy to obtain the result is reflected in the following schema.
\begin{center}
\begin{tabular}{ccc}
$C_\EA(\EA)$ & $\equiv$ & $\EA + \Rfn(\EA)$ \\
&&\\
$\rotatebox[origin=c]{90}{$\equiv$}_{\Pi_1}$ &  & $\rotatebox[origin=c]{90}{$\equiv$}_{\Pi_1}$ \\
&&\\
$\EA^+$ & \ \ \ $\equiv_{\Pi_1}$ & $\EA_\omega$ \\
\end{tabular}
\end{center}

The following two lemmas show that the vertical equivalences in the above schema are provable in $\EA$.
Moreover, they produce the functions effecting the proof-transformation.
\begin{lemma}\label{lm:poly}
$\EA \vdash \forall \pi \in \Pi_1\, (\Box_{\EA^+} \pi \eqv \Box_\EA [1]_\EA \pi)$.
\end{lemma}
\begin{proof}
Arguing in $\EA$ assume $\EA^+ \vdash \pi$ for some $\Pi_1$-sentence $\pi$. It follows that
$\EA \vdash \la 1\ra_\EA \top \imp \pi$, whence $\EA \vdash \la 1\ra_\EA \top \imp [1]_\EA \pi$ by provable $\Sigma_2$-completeness.
But clearly $\EA \vdash [1]_\EA \bot \imp [1]_\EA\pi$, so $\EA \vdash [1]_\EA \pi$.

Conversely, assume $\EA \vdash [1]_\EA \pi$ for some $\Pi_1$-sentence $\pi$. We derive
\begin{align*}
\EA \vdash \neg \pi &\imp [1]_\EA \neg \pi \wedge [1]_\EA \pi\\
&\imp [1]_\EA\bot,
\end{align*}
whence $\EA^+ \equiv  \EA + \la 1\ra_\EA \top  \vdash \pi$.
\end{proof}

For a proof of the next lemma we refer to \cite[Proposition 6.1]{Bek99b}.
\begin{lemma}
$\EA \vdash \forall \pi \in \Pi_1\, (\Box_{\EA + \Rfn(\EA)} \pi \eqv \Box_{\EA_{\omega}} \pi)$.
\end{lemma}

It follows from the proof of Lemma \ref{lm:poly} that there is a polynomially bounded function $f_1(x)$ that, given an $\EA^+$-proof $p$ of a $\Pi_1$-sentence $\pi$,
transforms it into an $\EA$-proof $f(p)$ of $[1]_\EA\pi$. Similarly, there is an elementary function $f_2(x)$ that acts the same
relative to the pair of theories $\EA + \Rfn(\EA)$ and $\EA_\omega$.
The following lemma shows that there is no such an elementary function for the pair of theories $\EA^+$ and $\EA_\omega$.

\begin{lemma} \label{lm:speedup}
$\EA^+$ has superexponential speed-up over $\EA_\omega$ w.r.t. $\Pi_1$-sentences.
\end{lemma}
\begin{proof}
We will show that $\EA^+$ proves the consistency of $\EA_\omega$ on some cut in $\EA^+$, whence the result follows by mimicking the proof
of \cite[Theorem 4.2]{Pud86}. Namely, one can take $\varphi_n := \Con_{\EA_\omega}(2^1_n)$,
where $\Con_{\EA_\omega}(x)$ is the formalization of ``there is no $\EA_\omega$-proof of contradiction of size $\leqslant x$''.

Define the formula $J(x) := \Con(\EA_x)$. We claim that $J(x)$ is the required cut in $\EA^+$.
Indeed, $\EA^+ \vdash J(0) \wedge \forall y \leqslant x\, (J(x) \imp J(y))$ is trivial.
Now we show that
$$
\EA^+ \vdash \forall x\, (J(x) \imp J(x + 1)).
$$
Using the fact that $\EA_{x+1} \equiv \EA + \Con(\EA_x)$ and $\EA^+ \equiv \EA + \la 1\ra_\EA\top$, we derive
\begin{align*}
\EA^+ \vdash \neg \Con(\EA_{x+1}) &\imp \Box_{\EA}(\neg \Con(\EA_{\overline{x}}))\\
&\imp \neg \Con(\EA_x),
\end{align*}
that is, $\EA^+ \vdash \neg J(x+1) \imp \neg J(x)$, whence $\EA^+ \vdash J(x) \imp J(x + 1)$.

Furthermore, we have
$$
\EA \vdash \forall x\, (J(x) \imp \Con_{\EA_\omega}(x)).
$$
Reasoning in $\EA$ assume that some $y \leqslant x$ codes an $\EA_\omega$-proof of a contradiction.
By the definition of $\EA_\omega$ the number $y$ must code a proof of contradiction in $\EA_{z+1}$ for some $z$
(namely, the maximal number $k$ such that the axiom $\Con(\EA_k)$ is used in the proof coded by $y$).
In particular, $\EA_{z+1}$ is inconsistent. By the definition of $z$ the proof coded by $y$ should contain the axiom $\Con(\EA_z)$,
and since $z \leqslant \ul \bar{z} \ur \leqslant \ul \Con(\EA_{\bar{z}})\ur$ we have $z < y \leqslant x$.
In this case the inconsistency of $\EA_{z+1}$ implies the inconsistency of $\EA_x$, whence $\neg \Con(\EA_x)$, that is, $\neg J(x)$.
\end{proof}

\begin{theorem}\label{th:speedup}
$C_\EA(\EA)$ has superexponential speed-up over $\EA + \Rfn(\EA)$.
\end{theorem}
\begin{proof}
By Lemma \ref{lm:speedup} there is a sequence of bounded sentences $\varphi_n$ with short $\EA^+$-proofs and long $\EA_\omega$-proofs
in a sense of the definition of speed-up.
We claim that this sequence also witnesses the superexponential speed-up of $C_\EA(\EA)$ over $\EA + \Rfn(\EA)$,
since the functions $f_1(x)$ and $f_2(x)$ preserve the bounds from the definition of speed-up.

Indeed, the sentences $\varphi_n$ have short  $C_\EA(\EA)$-proofs by transforming their short $\EA^+$-proofs via $f_1(x)$, which is polynomially bounded.
However, if there were short proofs of these sentences in $\EA + \Rfn(\EA)$ then we could transform them via the function $f_2(x)$, which is multiexponentially bounded,
into short $\EA_\omega$-proofs, contradicting Lemma \ref{lm:speedup}.
\end{proof}

We also prove the generalizations of the previous results.

\begin{lemma}\label{lm:rel-speedup}
For all $n \geqslant 0$ and $m > 0$ we have
$$
\EA \vdash \forall \pi \in \Pi_1\, (\Box_{\ISi_m} \pi \imp \Box_{\ISi_m} [n+1]_\EA \pi).
$$
\end{lemma}
\begin{proof}
We have $\EA \vdash \forall \pi \in \Pi_1\, (\Box_{\ISi_m}(\pi \imp [n+1]_\EA \pi))$ by provable $\Sigma_{n+1}$-completeness, whence the result follows.
\end{proof}

The following lemma is \cite[Proposition 6.2]{Bek99b}.
\begin{lemma}\label{lm:reflconserv}
For any $n \geqslant 0$ we have
$$
\EA \vdash \forall \alpha \succ 0\, \forall \pi \in \Pi_1\, (\Box_{\Rfn^n(\EA)_\alpha} \pi \eqv \Box_{(\EA)^n_{\omega^\alpha}} \pi).
$$
\end{lemma}

\begin{lemma} \label{lm:speedup2}
$\ISi_n$ has superexponential speed-up over $(\EA)^n_\omega$ w.r.t. $\Pi_1$-sentences.
\end{lemma}
\begin{proof}
Essentially the same as that of Lemma \ref{lm:speedup}, but with $J(x) := \nCon((\EA)^n_x)$ and using that $\ISi_n \equiv \EA + \la n+1 \ra_\EA\top$.
\end{proof}

\begin{theorem}
$C^{n+1}_{\EA}(\ISi_n)$ has superexponential speed-up over $\EA + \Rfn^n(\EA)$.
\end{theorem}
\begin{proof}
The proof is the same as that of Theorem \ref{th:speedup} by using the previous lemmas.
\end{proof}

\begin{corollary}
$C_{\EA}(\ISi_n)$ has superexponential speed-up over $\Rfn(\EA)_{\omega_n}$.
\end{corollary}
\begin{proof}
All the inclusions in the following chain are provable in $\EA$
$$
\ISi_n \subseteq_{\Pi_1} C_\EA(\ISi_n) \equiv \Rfn(\EA)_{\omega_n} \subseteq_{\Pi_1} \EA_{\omega_{n+1}} \subseteq
(\EA)^1_{\omega_n} \subseteq \dots \subseteq(\EA)^n_\omega.
$$
Indeed, for the inclusion $\ISi_n \subseteq_{\Pi_1} C_\EA(\ISi_n)$ it follows from Lemma \ref{lm:rel-speedup},
for $\Rfn(\EA)_{\omega_n} \subseteq_{\Pi_1} \EA_{\omega_{n+1}}$ it follows from Lemma \ref{lm:reflconserv}. As for the rest, it can be seen
from the proof of \cite[Theorem 3]{Bek99b} that the corresponding inclusions can be proved in $\EA$.

It follows that there are corresponding elementary, hence multiexponentially bounded, functions effecting the proof transformation
and for the first inclusion we obtain a polynomially bounded function from the proof of Lemma \ref{lm:rel-speedup}.
The rest of the proof is exactly the same as in Theorem~\ref{th:speedup}, but now using Lemma \ref{lm:speedup2}.
\end{proof}

\bibliographystyle{plain}
\bibliography{ref-all2}

\end{document}